\documentclass[12pt,a4paper]{article}

\usepackage{amsmath,amsfonts}
\usepackage{mathrsfs}

\newtheorem{proposition}{Proposition}

\let\scr\mathscr

\def\zs#1{_{\lower 3pt \hbox{$\scriptstyle#1$}}}
\def\Pb{\mathbf{P}}

\def\Ex{\mathbf{E}}

\def\sgn{{\rm sgn}} 

\begin{document}
\title{On the Goodness-of-Fit Testing for Ergodic Diffusion Processes}
\author{Yury A. \textsc{Kutoyants}\\
{\small Laboratoire de Statistique et Processus, Universit\'e du Maine}\\
{\small 72085 Le Mans, C\'edex 9, France}}

\date{}

\maketitle
\begin{abstract}
We consider the goodness of fit testing problem for ergodic diffusion
processes. The basic hypothesis is supposed to be  simple.  The diffusion
coefficient is known and the alternatives are described by the different trend
coefficients. We study the asymptotic distribution of the  Cramer-von Mises
type tests based on the empirical distribution function and local time
estimator of the invariant density. At particularly, we propose a
transformation which makes these tests asymptotically distribution free. We
discuss the modifications of this test in the case of composite basic
hypothesis. 

\end{abstract}


\section{Introduction}

The goodness of fit (GoF) tests play a special role in statistics because they
form a bridge between mathematical models and the real data. In classical
situation of i.i.d. observations $X^n=\left\{X_1,\ldots,X_n\right\}$ and the
basic hypothesis ${\scr H}_0$ : {\it the distribution function of $X_j$ is
$F_0\left(x\right)$}, the traditional solution is to construct a test
statistics $\Delta _n=D\left(\hat F_n,F_0\right)$ based on some distance
between the empirical distribution function $\hat F_n\left(x\right)$ and the
given (known) function $F_0\left(x\right)$. Then the test function is defined
by $\Psi _n=1_{\left\{\Delta _n >c_\varepsilon \right\}}$, where the constant
$c_\varepsilon $ is chosen from the condition: $\lim_{n\rightarrow \infty }
\Pb_0\left\{\Delta _n >c_\varepsilon\right\}=\varepsilon $, $\varepsilon \in
\left(0,1\right)$. The diversity of tests comes from the diversity of
distances. At particularly, if we take
\begin{equation}
\label{1}
\Delta_n =n\int_{-\infty }^{\infty }H\left(x\right)\left[\hat
F_n\left(x\right)-F_0\left(x\right)\right]^2{\rm d}F_0\left(x \right),
 \end{equation}
then we obtain the well-known  Cram\'er-von Mises family of statistics
\cite{LR}. If $H\left(x\right)\equiv 1$, then we have Cram\'er-von Mises test
and if the weight function   $H\left(x\right)=
\left(F_0\left(x\right)\left[1-F_0\left(x\right)\right]\right)^{-1} $ we
obtain  the  Anderson-Darling test. 

In the case of uniform metric 
\begin{equation}
\label{ks}
\Delta  _n=\sup_x\sqrt{n}\left| \hat
F_n\left(x\right)-F_0\left(x\right)\right|, 
\end{equation}
and  $\Psi _n=1_{\left\{\Delta _n >c_\varepsilon \right\}}$ we have
Kolmogorov-Smirnov test.  Remind that the tests based on these statistics are
asymptotically distribution free (ADF) and for continuous $F_0\left(x\right)$ under
hypothesis ${\scr H}_0$ we have the convergence
\begin{align*}
&n\int_{-\infty }^{\infty }\left[\hat
F_n\left(x\right)-F_0\left(x\right)\right]^2{\rm d}F_0\left(x
\right)\Longrightarrow \int_{0}^{1}W_0\left(s\right)^2\;{\rm d}s,\\
&n\int_{-\infty }^{\infty }\frac{\left[\hat
F_n\left(x\right)-F_0\left(x\right)\right]^2}{F_0\left(x\right)\left[1-F_0\left(x\right)\right]  }\;{\rm d}F_0\left(x
\right)\Longrightarrow
\int_{0}^{1}\frac{W_0\left(s\right)^2}{s\left(1-s\right)}\;{\rm d}s,\\ 
&\sup_x\sqrt{n}\left| \hat
F_n\left(x\right)-F_0\left(x\right)\right|\Longrightarrow \sup_{0\leq s\leq
1}\left|W_0\left(s\right) \right| ,
\end{align*}
where $W_0\left(x\right)$ is Brownian bridge. These last property of the
statistics allows us to chose once the constant $c_\varepsilon $ for all
$F_0\left(x\right)$. Note as well that the both statistics tend to $\infty $
for any fixed alternative $F\left(\cdot \right)\not=F_0\left(\cdot \right)$
and this property provides the consistency of these tests. 

The present work is devoted to the similar problem but in the case of
continuous time observations $X^T=\left\{X_t,0\leq t\leq T\right\}$ of ergodic diffusion process
\begin{equation}
\label{2}
{\rm d}X_t=S\left(X_t\right)\;{\rm d}t +\sigma
\left(X_t\right)\;{\rm d}W_t,\qquad X_0, \quad 
0\leq t\leq T .
\end{equation}
The diffusion coefficient $\sigma \left(\cdot \right)^2$ is supposed to be
known and the hypothesis is concern the trend coefficient $S\left(\cdot
\right)$ only. That means, that the basic hypothesis is simple:\\ ${\scr H}_0$
: {\it The observed trajectory $X^T$ is
solution of the stochastic differential equation}
$$
{\rm d}X_t=S_0\left(X_t\right)\;{\rm d}t +\sigma \left(X_t\right)\;{\rm
d}W_t,\qquad X_0, \quad 0\leq t\leq T 
$$  
{\it where $S_0\left(x\right)$ is some known function}.

The alternative corresponds
to the case $S\left(\cdot \right)\not = S_0\left(\cdot \right)$. 

We suppose that the trend (under hypothesis and alternative) and diffusion
coefficients satisfy the conditions:

${\cal ES}$. {\it The function $S\left(\cdot \right)$ is locally bounded, the
function $\sigma \left(\cdot \right)^2$ is continuous and positive and for
some $A>0$ the inequality
$$
xS\left(x\right)+\sigma \left(x\right)^2\leq A\left(1+x^2\right)
$$
holds.}

 This condition provides the existence of unique weak  solution of this
 equation (see \cite{DuR}, p. 210).

Moreover, we suppose that the  following condition is fulfilled too.

${\cal RP}$. {\it The function
$$
V\left(S,x\right)=\int_{0}^{x}\exp\left\{-2\int_{0}^{y}\frac{S\left(z\right)}{\sigma
\left(z\right)^2}\;{\rm d}z\right\}\;{\rm d}y\longrightarrow \pm \infty 
$$
as $x\rightarrow \pm \infty $ and}
$$
G\left(S\right)=\int_{-\infty }^{\infty }\sigma
\left(y\right)^{-2}\exp\left\{2\int_{0}^{y}\frac{S\left(z\right)}{\sigma
\left(z\right)^2}\;{\rm d}z\right\}\;{\rm d}y<\infty
$$

By this condition the observed process is recurrent positive and has ergodic
properties with the density of invariant law
$$
f\left(x\right)=\frac{1}{G\left(S\right)\,\sigma \left(x\right)^{2}
}\exp\left\{2\int_{0}^{x}\frac{S\left(y\right)}{\sigma 
\left(y\right)^2}\;{\rm d}y\right\}.
$$
The corresponding density and distribution function under hypothesis ${\scr
H}_0$ we denote as $f_0\left(x\right)$ and $F_0\left(x\right)$ and the
mathematical expectation as $\Ex_0$. Moreover, we suppose that the initial
value $X_0$ is a random variable with this distribution function because this
condition simplifies exposition (the observed process is stationary).

Let us fix some $\varepsilon \in \left(0,1\right)$ and denote by ${\cal
K}_\varepsilon $ the class of tests $\psi _T$ of asymptotic size  $\varepsilon
$, i.e.; $\Ex_0\psi _T=\varepsilon +o\left(1\right)$. We are interested by the GoF tests of asymptotic size $\varepsilon $, which
are ADF. 

The problem of goodness of fit testing can be considered as follows: let us
introduce some statistic $\delta _T$ such that its limit distribution
$G\left(x\right)$ under hypothesis does not depend on the model, then the test
$\psi _T=1_{\left\{\delta _T>c_\varepsilon \right\}}$, where  the constant
$c_\varepsilon $ is solution of the equation $1-G\left(c_\varepsilon
\right)=\varepsilon $ is ADF. Moreover, we
require as well that for any fixed alternative (defined by the trend
coefficient $S\left(\cdot \right)$) we have $\Pb\left\{\delta _T>c_\varepsilon
\right\}\rightarrow 1$, i.e.; the test is consistent.  

Let us remind here some of  known tests satisfying these conditions. 
 The first ones
$
\psi_T\left(X^T\right)=1_{\left\{\delta _T>d_\varepsilon \right\}}$ and $
\phi_T\left(X^T\right)=1_{\left\{\gamma _T>c_\varepsilon\right\}} 
$
 are based on the
   following two statistics
\begin{align*}
 \delta _T&=\frac{1}{T^2\,\Ex_{0}\left[\sigma \left(\xi \right)^2\right]}\;
\int_{0}^{T}\left[X_t-X_0-\int_{0}^{t} S_{0}\left(X_v\right)\;{\rm d}v
\right]^2{\rm d}t,\\
\gamma _T&=\frac{1}{\sqrt{T\,\Ex_{0}\left[\sigma \left(\xi \right)^2\right]}}\;
\sup_{0\leq t\leq T}\left|X_t-X_0-\int_{0}^{t} S_0\left(X_v\right)\;{\rm d}v\right|.
\end{align*}
Here and in the sequel $\xi $ is the random variable with the density
$f_0\left(\cdot \right)$. 

It is shown that under hypothesis ${\scr H}_0$
\begin{align*}
\delta_T\Longrightarrow \int_{0}^{1}w\left(v\right)^2{\rm d}v,\qquad 
\gamma _T\Longrightarrow  \sup_{0\leq v\leq 1}\left|w\left(v\right)\right|.
\end{align*}
Hence the  constants $c_\varepsilon $ and $d_\varepsilon $ are defined by the
equations
\begin{align}
\label{4a}
\Pb\left\{\int_{0}^{1}w\left(v\right)^2{\rm d}v>d_\varepsilon
\right\}=\varepsilon ,\qquad \Pb\left\{\sup_{0\leq v\leq
1}\left|w\left(v\right)\right|>c_\varepsilon \right\}=\varepsilon .
\end{align}.
  These tests belong to ${\cal K}_\varepsilon $ and are consistent
against any alternative  $S\left(\cdot \right)\not=S_0\left(\cdot \right)$
such that $\Ex_{0}S\left(\xi\right)  \not = 0 $.
Moreover, it is shown that the test  $\psi_T\left(X^T\right)$ is
asymptotically optimal in special sense (see \cite{DK} for details).

Another GoF test $\varphi _T=1_{\left\{\hat\Delta _T>e_\varepsilon \right\}}$ was
proposed by Negri and Nishiyama \cite{NeNi}. It is based on the statistic
\begin{equation}
\label{nn}
\hat\Delta_ T=\frac{1}{\sqrt{T\Ex_{0}\sigma \left(\xi
\right)^2   }}\;\ \sup_x\left|\int_{0}^{T}1_{\left\{X_t<x\right\}}\;\left[{\rm 
d}X_t-S_0\left(X_t\right){\rm d}t\right] \right|
\end{equation}
which converges to $\hat\Delta_0=\sup_{0\leq v\leq 1}
\left|w\left(v\right)\right|$. The constant $c_\varepsilon $ is defined in
\eqref{4a}.  This test belongs to ${\cal K}_\varepsilon $ and is consistent
against any fixed alternative satisfying condition: for some $x$ we have $
\Ex_{0}\left(1_{\left\{\xi <x\right\}}\left[S\left(\xi \right)-S_0\left(\xi
\right)\right]\right) \not = 0.  $
 
Note that the similar question of the goodness of fit testing for ergodic
diffusion processes by discrete time observations was extensively studied (see
Chen and Gao \cite{CG} and references therein). 

The goal of this work is to study the tests which are direct analogues  of the
classical GoF tests like Anderson-Darling \eqref{1} and Kolmogorov-Smirnov
\eqref{ks}.

To test the hypothesis ${\scr H}_0$ we propose two tests of Cram\'er-von Mises
type. The first one is based on {\it empirical distribution function} (EDF)
$$
\hat F_T\left(x\right)=\frac{1}{T}\int_{0}^{T}1_{\left\{X_t<x\right\}}\,{\rm d}t 
$$
and the statistic is similar to \eqref{1}:
$$
\Delta _T=T\int_{-\infty }^{\infty }H\left(x\right)\left[\hat
F_T\left(x\right)-F_0\left(x\right)\right]^2 {\rm d}F_0\left(x\right).
$$
The second test is based on {\it local time estimator} (LTE) $\hat
f_T\left(x\right)$ of the invariant density, which can be written as
$$
\hat f_T\left(x\right)=\frac{\Lambda _T\left(x\right)}{T\,\sigma
\left(x\right)^2}= \frac{\left|X_T-x\right|-\left|X_0-x\right|}{T\,\sigma
\left(x\right)^2 }-\frac{1}{T\,\sigma
\left(x\right)^2}\int_{0}^{T}\sgn\left(X_t-x\right)\;{\rm d}X_t. 
$$
Here $\Lambda _T\left(x\right)$ is the local time of the diffusion
process. The corresponding statistic is
$$
\delta _T=T\int_{-\infty }^{\infty }h\left(x\right)\left[\hat
f_T\left(x\right)-f_0\left(x\right)\right]^2 {\rm d}F_0\left(x\right).
$$
We discuss as well the Kolmogorov-Smirnov type test with the test statistic
\begin{equation}
\label{Ks}
\gamma _T=\sup_x \sqrt{T}g\left(x\right)\left|\hat
f_T\left(x\right)-f_0\left(x\right) \right|.
\end{equation}
The goodness of fit tests are $\Psi _T=1_{\left\{\Delta _T>c_\varepsilon \right\}}$,
and $\psi _T=1_{\left\{\delta _T>d_\varepsilon \right\}}$ and $\hat\psi
_T=1_{\left\{\gamma  _T>e_\varepsilon \right\}}$. These tests with 
$H\left(x\right)\equiv 1, h\left(x\right)\equiv 1,g\left(x\right)\equiv 1$, were 
proposed  in \cite{DK}, but  unfortunately they  are not ADF (see
as well \cite{Fou}, \cite{FouK}, where  similar test statistics  are
discussed). 

The test $\varphi _T$, \eqref{nn} is ADF. It is interesting to see the
relation between the statistics \eqref{nn} and \eqref{Ks}. Suppose that
$\sigma \left(x\right)\equiv 1$ and remind that
$$
\bar f_T\left(x\right)=\frac{2}{T}\int_{0}^{T}1_{\left\{X_t<x\right\}}\;{\rm
d}X_t,\qquad \quad
f_T^*\left(x\right)=\frac{2}{T}\int_{0}^{T}1_{\left\{X_t<x\right\}}
S_0{\left(X_t\right)}\;{\rm d}t 
$$ 
are unbiased estimators of the invariant density \cite{Kut97}. The first is
asymptotically equivalent to the local time estimator $\hat
f_T\left(x\right)$, i.e.; $$\sqrt{T}\left(\bar f_T\left(x\right)- \hat
f_T\left(x\right)\right)=o\left(1\right) $$ and for the second we have
\begin{align*}
f_T^*\left(x\right)\longrightarrow 2\Ex_01_{\left\{\xi
<x\right\}}S_0\left(\xi  \right) =f_0\left(x\right).
\end{align*}
Of course, the second estimator is not indeed estimator of the invariant
density because it uses the function $S_0\left(x\right)$, but this choice of
the second term makes the statistic asymptotically distribution free. 
Therefore
$$
\hat \Delta _T=\frac{\sqrt{T}}{2}\sup_x\left|\bar
f_T\left(x\right)-f_T^*\left(x\right)\right| .  
$$

  The random processes $\left\{\eta _T\left(x\right), x\in {\cal R}\right\}$
  and $\left\{\zeta _T\left(x\right), x\in {\cal R}\right\}$ where
$$
\eta _T\left(x\right)=\sqrt{T}\left[\hat
F_T\left(x\right)-F_0\left(x\right)\right]\quad {\rm  and}\quad 
\zeta _T\left(x\right)=\sqrt{T}\left[\hat
f_T\left(x\right)-f_0\left(x\right)\right]
$$ 
converge to
the Gaussian processes $\left\{\eta \left(x\right), x\in {\cal R}\right\}$ and
$\left\{\zeta \left(x\right), x\in {\cal R}\right\}$ with quite complicate
structure (for the first convergence see \cite{Neg} and for the second see
\cite{Kut04}). For example,
$$
\Ex_0\zeta \left(x\right)\zeta
\left(y\right)=4f_0\left(x\right)f_0\left(y\right)\int_{-\infty }^{\infty } \frac{\left(
1_{\left\{z>x\right\}}-F_0\left(z\right)\right) \left(
1_{\left\{z>y\right\}}-F_0\left(z\right)\right)  }{\sigma
\left(z\right)^2f_0\left(z\right)}\;  {\rm d}z
$$

The goal of this work is to show that for certain choice of weight functions
$H\left(x\right)$ and $h\left(x\right)$ the tests $\Psi _T=1_{\left\{\Delta
_T>c_\varepsilon \right\}}$ and $\psi _T=1_{\left\{\delta _T>d_\varepsilon \right\}}$
can be ADF and consistent  against some alternatives.

We discuss as well some other  ADF GoF tests.

\section{GoF tests based on LTE}

Let us denote by $\mu $ the median of the invariant law : $F_0\left(\mu
\right)=1/2$.
We consider a special hypotheses testing problem, when the trend coefficient
is fixed for $x\leq \mu $ or for $x\geq \mu $. Therefore, we consider {\it one
sided alternatives}  only. Suppose that the changes in the trend coefficient
are possible for $x\geq \mu $, i.e., the values of $S\left(x\right)$ for
$x<\mu $ are the same under hypothesis and alternative. 

We study first the Cram\'er-von Mises type test $\psi _T=1_{\left\{\delta
_T>d_\varepsilon \right\}}$, where the statistic
$$
\delta_T=T\int_{\mu  }^{\infty } h\left(x\right) \left(\hat f _T
\left(x\right) -f_0\left(x\right)\right)^2{\rm d}F_0\left(x \right)
$$
with
$$
h\left(x\right)=\frac{2F_0\left(x\right)-1 }{4\Phi \left(\mu \right)^2\sigma
\left(x\right)^2f_0\left(x\right)^4}\,\;e^{-\Phi \left(x\right)/\Phi
\left(\mu \right)}\;1_{\left\{x\geq \mu \right\}},
$$
and
$$
\Phi \left(x\right)=\int_{-\infty }^{\infty } \frac{\left(
1_{\left\{y>x\right\}}-F_0\left(y\right)\right)^2}{\sigma
\left(y\right)^2f_0\left(y\right)}\;  {\rm d}y.
$$
The study of this test is based on the asymptotic normality of the LTE:
$$
\sqrt{T}\left(\hat f_T\left(x\right)- 
f_0\left(x\right)\right)\Longrightarrow \zeta \left(x\right)\sim  {\cal N}\left(0,d_f\left(x\right)^2\right),
$$
where $d_f\left(x\right)^2=4f_0\left(x\right)^2\Phi \left(x\right)$.
But we have two other classes of density estimators, 
 which are consistent and asymptotically normal with the same limit
distribution. These are kernel type estimators 
$$
\bar
f_T\left(x\right)=\frac{1}{\sqrt{T}}\int_{0}^{T}
K\left({\sqrt{T}\left(X_t-x\right)}\right)\;{\rm d}t ,
$$
and the unbiased estimators
$$
\tilde{f}_T(x)=\frac{2}{ T }\int _0^T \frac{1_{\{X_t\leq x\}}\; h(X_t)}{
\sigma (x)^2\, h(x)  }\;{\rm d}X_t+\frac{1}{ T }\int _0^T \frac{1_{\{X_t\leq
x\}}\; h'(X_t)\;\sigma (X_t)^2}{ \sigma (x)^2\, h(x) }\;{\rm d}t
$$
where $h\left(\cdot \right)\in {\cal C}'$ is an arbitrary function. Under mild
regularity conditions (see \cite{Kut04}, Propositions 1.58 and 1.61) we have
the same asymptotic normality
$$
\sqrt{T}\left(\bar f_T\left(x\right)- 
f\left(x\right)\right)\Longrightarrow  {\cal
N}\left(0,d_f\left(x\right)^2\right)
$$
and
$$
 \sqrt{T}\left(\tilde
f_T\left(x\right)-  
f\left(x\right)\right)\Longrightarrow  {\cal N}\left(0,d_f\left(x\right)^2\right).
$$
Therefore, the results obtained for the test $\psi  _T$ based on the LTE can
be valid  for the tests based on the kernel-type or unbiased estimators as
well (under strengthen conditions). 

The calculation of the statistic $\delta _T$ can be slightly
simplified. Remind the equalities
$$
\frac{1}{T}\int_{0}^{T}g\left(X_t\right)\,{\rm d}t=\int_{-\infty }^{\infty
}g\left(x\right)\, \hat f_T\left(x\right)\,{\rm d}x =\int_{-\infty }^{\infty
}g\left(x\right)\, \,{\rm d}\hat F_T\left(x\right).
$$
Hence
\begin{align*}
\delta _T&=T\int_{\mu  }^{\infty } h\left(x\right) \left(\hat f _T
\left(x\right) -f_0\left(x\right)\right)^2{\rm d}\hat F_T\left(x \right)\\
&\quad +T\int_{\mu  }^{\infty } h\left(x\right) \left(\hat f _T
\left(x\right) -f_0\left(x\right)\right)^2{\rm d}\left[F_0\left(x\right)- \hat
F_T\left(x \right)\right]\\
&= \int_{0}^{T}h\left(X_t\right)\left[\hat f _T
\left(X_t\right) -f_0\left(X_t\right) \right]^2   \,{\rm
d}t-\frac{1}{\sqrt{T}} \int_{\mu  }^{\infty }h\left(x\right)\,\zeta
_T\left(x\right)^3{\rm d}x \\
&= \int_{0}^{T}h\left(X_t\right)\left[\hat f _T
\left(X_t\right) -f_0\left(X_t\right) \right]^2   \,{\rm
d}t+o\left(1\right)
\end{align*}
and we can use the statistic
\begin{equation}
\label{emp}
\delta _T^*=\int_{0}^{T}h\left(X_t\right)\left[\hat f _T \left(X_t\right)
-f_0\left(X_t\right) \right]^2 \,{\rm d}t
\end{equation}
in the construction of our GoF test. To avoid the difficulties with
calculation of stochastic integral related to $\hat f _T \left(X_t\right) $ we
can use here the kernel type estimator of the density.

Note that the LTE of the density is the derivative with probability 1 of the
EDF. Indeed, we have the equality ($r>0$)
$$
\frac{\hat F_T\left(x+r\right)-\hat
F_T\left(x\right)}{r}=\frac{1}{rT}\int_{0}^{T}1_{\left\{ x\leq X_t\leq
x+r\right\}}\;{\rm d}t =\frac{1}{r}\int_{x}^{x+r}\hat f_T\left(y\right)\;{\rm
d}y
$$
and, as the local time is continuous with probability one (see \cite{RY}), we have
$$
\lim_{r\rightarrow 0}\frac{1}{r}\int_{x}^{x+r}\hat f_T\left(y\right)\;{\rm
d}y =\hat f_T\left(x\right).
$$

\bigskip

Let us define the constant $d_\varepsilon $ as solution of the equation
$$
\Pb\left\{\int_{1}^{\infty
}w{\left(v\right)}^2\;e^{-v}\;{\rm
d}v>d_\varepsilon  \right\}=\varepsilon,
$$
where $w\left(\cdot \right)$ is a Wiener process and introduce the following condition.

{\bf Condition} ${\bf \cal A}$. {\it The function $\Phi \left(x\right)<\infty $ for all
$x\in \left(\mu ,\infty \right)$ and} 
\begin{align}
\label{3}
A_1&= \int_{\mu }^{\infty }\frac{\left(2F_0\left(x\right)-1\right)\Phi
\left(x\right)}{\Phi \left(\mu \right)^2\sigma 
\left(x\right)^2f_0\left(x\right)^2}\;  \,e^{-\frac{\Phi
\left(x\right)}{\Phi \left(\mu \right)}}\;\,{\rm d}x<\infty,\\
\label{4}
A_2&= \int_{\mu }^{\infty }\frac{\left(2F_0\left(x\right)-1\right)\,e^{-\frac{\Phi
\left(x\right)}{\Phi \left(\mu \right)}}}{
\Phi \left(\mu \right)^2\sigma 
\left(x\right)^2f_0\left(x\right)^2}\; \Ex_0 \left(\int_{0}^{\xi }\frac{
1_{\left\{v>x\right\}}- F_0\left(v\right)}{\sigma
\left(v\right)^2f_0\left(v\right)}\;{\rm d}v\right)^2 {\rm d}x<\infty.
\end{align}

\begin{proposition}
\label{P1} 
Let the conditions ${\cal ES}$, ${\cal RP}$ and ${\cal A}$ be fulfilled, then
the GoF test $\psi _T=1_{\left\{\delta _T>d_\varepsilon \right\}}$ is
ADF.
\end{proposition}

{\bf Proof.}
We have to prove the convergence (under hypothesis ${\scr H}_0$)
\begin{align*}
\delta_T\Longrightarrow \delta _0&=\int_{1}^{\infty }w\left(v\right)^2e^{-v}\,{\rm d}v .
\end{align*}

First we remind some properties of the local time estimator of the density
(under hypothesis ${\scr H}_0$). 
The following presentation
is valid (see \cite{Kut04},   p. 29)
\begin{align}
\zeta _T\left(x\right)&=
\frac{2f_0\left(x\right)}{\sqrt{T}}\int_{0}^{T}\frac{ F_0\left(X_t\right)-
1_{\left\{X_t>x\right\}}}{\sigma \left(X_t\right)f_0\left(X_t\right)}\;{\rm
d}W_t \nonumber\\ 
&\label{5}\qquad
+\frac{2f_0\left(x\right)}{\sqrt{T}}\int_{X_0}^{X_T}\frac{
1_{\left\{v>x\right\}}-  F_0\left(v\right)}{\sigma
\left(v\right)^2f_0\left(v\right)}\;{\rm d}v=\zeta
_T^{\left(1\right)}\left(x\right)+\zeta _T^{\left(2\right)}\left(x\right)  
\end{align}
in obvious notation. 
We supposed that the process $X_t,t\geq 0$ is stationary, hence
\begin{align*}
\Ex_0 \left(\zeta _T^{\left(2\right)}\left(x\right)\right)^2&=
\Ex_0\left(\frac{2f_0\left(x\right)}{\sqrt{T}}\int_{X_0}^{X_T}\frac{ 
1_{\left\{v>x\right\}}-  F_0\left(v\right)}{\sigma
\left(v\right)^2f_0\left(v\right)}\;{\rm d}v\right)^2\\
& \leq \frac{16}{T}\;\Ex_0
\left(\int_{0}^{\xi }\frac{ 
1_{\left\{v>x\right\}}- F_0\left(v\right)}{\sigma
\left(v\right)^2f_0\left(v\right)}\;{\rm d}v\right)^2 .
\end{align*}
Further,  by condition \eqref{4} we have
\begin{align*}
&4\int_{\mu }^{\infty }h\left(x\right)f_0\left(x\right)^2\Ex_0
\left(\int_{0}^{\xi }\frac{ 
1_{\left\{v>x\right\}}- F_0\left(v\right)}{\sigma
\left(v\right)^2f_0\left(v\right)}\;{\rm d}v\right)^2{\rm d}x\\
&= \int_{\mu }^{\infty }\frac{\left(2F_0\left(x\right)-1\right)\,e^{-\frac{\Phi
\left(x\right)}{\Phi \left(\mu \right)}}}{
\Phi \left(\mu \right)^2\sigma 
\left(x\right)^2f_0\left(x\right)^2}\; \Ex_0 \left(\int_{0}^{\xi }\frac{
1_{\left\{v>x\right\}}- F_0\left(v\right)}{\sigma
\left(v\right)^2f_0\left(v\right)}\;{\rm d}v\right)^2 {\rm d}x<\infty .
\end{align*}
Hence
\begin{align*}
 \int_{\mu }^{\infty }h\left(x\right)\Ex_0\left(\zeta _T^{(2)}\right)^2{\rm
 d}F_0\left(x\right)\leq \frac{A_2}{T}\longrightarrow 0.
\end{align*}
Therefore it is sufficient to study the first integral in \eqref{5}, i.e.,
$$
\delta _T=\int_{\mu }^{\infty
}\frac{\left(2F_0\left(x\right)-1\right)\,e^{-\frac{\Phi \left(x\right)}{\Phi
\left(\mu \right)}}}{T \Phi \left(\mu \right)^2\sigma
\left(x\right)^2f_0\left(x\right)^2}\;\left(\int_{0}^{T}\frac{1_{\left\{X_t>x\right\}}-
F_0\left(X_t\right) }{\sigma \left(X_t\right)f_0\left(X_t\right)}\;{\rm d}W_t
\right)^2{\rm d}x+o\left(1\right) .
$$
Let us denote
$$
\hat\zeta _T\left(x\right)=\frac{1}{\sqrt{T}}\int_{0}^{T}\frac{1_{\left\{X_t>x\right\}}-
F_0\left(X_t\right) }{\sigma \left(X_t\right)f_0\left(X_t\right)}\;{\rm d}W_t.
$$
By the law of large numbers we have the convergence
$$
\frac{1}{T}\int_{0}^{T}\left(\frac{
1_{\left\{X_t>x\right\}}-F_0\left(X_t\right)}{\sigma
\left(X_t\right){{f_0\left(X_t\right)}}}\right)^2\; {\rm
d}t\longrightarrow \Phi \left(x\right). 
$$
hence we can apply  the central limit theorem (see, e.g., \cite{Kut04},Theorem
1.19) and to obtain  the convergence of all finite dimensional distributions of $\hat\zeta
_T\left(\cdot \right)$ to the multidimensional Gaussian law:
\begin{equation}
\label{6}
\left(\hat\zeta _T\left(x_1\right),\ldots,\hat\zeta
_T\left(x_k\right)\right)\Longrightarrow \left(\hat\zeta
\left(x_1\right),\ldots,\hat\zeta\left(x_k\right)\right) ,
\end{equation}
where the Gaussian process $\zeta \left(x\right)$  is with zero mean and it's covariance
function is
$$
R\left(x,y\right)=\int_{-\infty }^{\infty } \frac{
\left(1_{\left\{v>x\right\}}-F_0\left(v\right)\right)
\left(1_{\left\{v>y\right\}}-F_0\left(v\right)\right)}{\sigma
\left(v\right)^2{{f_0\left(v\right)}}}\; {\rm d}v.
$$
The  process $\hat\zeta \left(x\right)$  admits the representation
$$
\hat\zeta \left(x\right) =\int_{-\infty }^{\infty } \frac{
1_{\left\{y>x\right\}}-F_0\left(y\right)}{\sigma
\left(y\right){\sqrt{f_0\left(y\right)}}}\; {\rm d}W_y.
$$
The last integral is with respect to 
double-sided Wiener process, i.e.; $W_y=W^+\left(y\right), y\geq
0$ and $W_y=W^-\left(-y\right),y\leq 0$, where $W^-\left(\cdot \right)$ and
$W^+\left(\cdot \right)$ are two independent Wiener processes. 

The convergence \eqref{6} formally can be obtained as follows (egalities in distribution)
\begin{align*}
\hat\zeta  _T\left(x\right)&=W\left(\frac{1}{T}\int_{0}^{T} \left(
\frac{\
1_{\left\{X_t>x\right\}}-F_0\left(X_t\right)}{\sigma
\left(X_t\right){f_0\left(X_t\right)}}\; \right)^2{\rm d}t \right)+o\left(1\right)\\
&=W\left(\int_{-\infty }^{\infty }  \left(
\frac{
1_{\left\{y>x\right\}}-F_0\left(y\right)}{\sigma
\left(y\right){f_0\left(y\right)}}\; \right)^2\hat f_T\left(y\right) {\rm d}t
\right)+o\left(1\right)\\
&\Longrightarrow  W\left(\int_{-\infty }^{\infty } \left(\frac{
1_{\left\{y>x\right\}}-F_0\left(y\right)}{\sigma
\left(y\right){f_0\left(y\right)}}\right)^2\; f_0\left(y\right) {\rm
d}y\right) \\
&=\int_{-\infty }^{\infty } \frac{
1_{\left\{y>x\right\}}-F_0\left(y\right)}{\sigma
\left(y\right){\sqrt{f_0\left(y\right)}}}\;  {\rm
d}W_y=\hat\zeta \left(x\right) 
\end{align*}
where $W\left(s \right),s\geq 0$ is some Wiener process.
 We used above the
following property of local time of ergodic diffusion process : for any
integrable function $g\left(\cdot \right)$
\begin{equation}
\label{7}
\frac{1}{T}\int_{0}^{T}g\left(X_t\right)\,{\rm d}t=\int_{-\infty }^{\infty
}g\left(y\right)\hat f_T\left(y\right)\, {\rm d} y\longrightarrow
\int_{-\infty }^{\infty 
}g\left(y\right)\; f_0\left(y\right)\, {\rm d} y .
\end{equation}
To verify the convergence
\begin{align*}
\delta_T\Longrightarrow 4\int_{\mu }^{\infty }
h\left(x\right)f_0\left(x\right)^2\hat\zeta \left(x\right)^2{\rm d}F_0\left(x
\right) =4\int_{\mu }^{\infty } h\left(x\right)f_0\left(x\right)^3W\left(\Phi
\left(x\right)\right)^2 {\rm d}x
\end{align*}
we note that for any $\varepsilon >0$ there exists $L>\mu $ such that
$$
\int_{L}^{\infty } h\left(x\right)\Ex_0\left(\zeta
_T^{\left(1\right)}\left(x\right)\right)^2{\rm d}F_0\left(x\right)\leq \varepsilon . 
$$
Therefore it is sufficient to show that for any $L>\mu $
$$
\int_{\mu }^{L } h\left(x\right)f_0\left(x\right)^2\,\hat\zeta
_T\left(x\right)^2{\rm
d}F_0\left(x\right)\Longrightarrow \int_{\mu }^{L }
h\left(x\right)f_0\left(x\right)^2\hat\zeta 
\left(x\right)^2{\rm d}F_0\left(x\right).
$$
This last convergence follows from the  the convergence of finite
dimensional distributions \eqref{6} and the estimate
\begin{equation}
\label{8}
\Ex_0\left|\hat\zeta _T\left(x_2\right)^2-\hat\zeta
_T\left(x_1\right)^2\right| \leq C_L\,\left|x_2-x_1\right|^{1/2}
\end{equation}
which is valid for all $\left|x_i\right|\leq L$ (see \cite{GS}, Theorem
9.7.1). The estimate \eqref{8} can be obtained as follows. We have
\begin{align*}
&\left(\Ex_0\left|\hat\zeta _T\left(x_2\right)^2-\hat\zeta
_T\left(x_1\right)^2\right| \right)^2\\
&\qquad \leq \Ex_0\left|\hat\zeta _T\left(x_2\right)-\hat\zeta
_T\left(x_1\right)\right|^2\Ex_0\left|\hat\zeta _T\left(x_2\right)+\hat\zeta
_T\left(x_1\right)\right|^2\\
&\qquad \leq \Ex_0\left|\hat\zeta _T\left(x_2\right)-\hat\zeta
_T\left(x_1\right)\right|^2\left(2\Ex_0\hat\zeta
_T\left(x_2\right)^2+2\Ex_0\hat\zeta
_T\left(x_1\right)^2\right)\\
&\qquad \leq \left(2\Phi \left(x_2\right)+2\Phi \left(x_1\right)  \right)\Ex_0\left|\hat\zeta _T\left(x_2\right)-\hat\zeta
_T\left(x_1\right)\right|^2.
\end{align*}
Further (let  $x_2>x_1$)
\begin{align*}
\Ex_0\left|\hat\zeta _T\left(x_2\right)-\hat\zeta
_T\left(x_1\right)\right|^2&=\Ex_0
\left(\frac{1}{\sqrt{T}}\int_{0}^{T}\frac{1_{\left\{x_1\leq X_t\leq
x_2\right\}}}{ \sigma \left(X_t\right)f_0\left(X_t\right)}\;{\rm
d}W_t\right)^2\\
 &\quad =\int_{x_1}^{x_2}\frac{{\rm d}v}{\sigma
\left(v\right)^2f_0\left(v\right)}\leq B_L\left|x_2-x_1\right|.
\end{align*}
Remind that $h\left(\cdot \right)$ is continuous on $\left[\mu
,L\right]$.

The last step is to verify the equality (in distribution)
$$
4\int_{\mu }^{\infty }
h\left(x\right)f_0\left(x\right)^2\hat\zeta \left(x\right)^2{\rm d}F_0\left(x
\right) =\int_{1}^{\infty }w\left(v\right)^2e^{-v}\,{\rm d}v.
$$

 For the function 
\begin{align*}
\Phi \left(x\right)&=\int_{-\infty }^{\infty } \frac{\left(
1_{\left\{y>x\right\}}-F_0\left(y\right)\right)^2}{\sigma
\left(y\right)^2f_0\left(y\right)}\;  {\rm d}y\\
&=\int_{-\infty }^{x } \frac{F_0\left(y\right)^2}{\sigma
\left(y\right)^2f_0\left(y\right)}\;  {\rm d}y+\int_{x}^{\infty } \frac{\left(
1-F_0\left(y\right)\right)^2}{\sigma
\left(y\right)^2f_0\left(y\right)}\;  {\rm d}y
\end{align*}
 we have
\begin{align*}
\Phi' \left(x\right)=\frac{2F_0\left(x\right)-1 }{\sigma
\left(x\right)^2f_0\left(x\right)}<0,\quad {\rm for}\quad x<\mu 
\end{align*}
and $\Phi' \left(x\right)>0$ for $x>\mu $. Hence the functions $\Phi
\left(x\right), x\leq \mu $ and $\Phi
\left(x\right), x\geq  \mu$ are strictly monotone (decreasing and increasing
respectively). Moreover, $\Phi \left(\pm\infty \right)=\infty $. 

We can write
\begin{align*}
&4\int_{\mu  }^{\infty } 
h\left(x\right)f_0\left(x\right)^3W\left(\Phi \left(x\right)\right)^2 {\rm
d}x\\
&\qquad=4\int_{\mu  }^{\infty }  
\frac{h\left(x\right)\sigma
\left(x\right)^2f_0\left(x\right)^4}{2F_0\left(x\right)-1}\;W\left(\Phi
\left(x\right)\right)^2\; \Phi '\left(x\right)\;{\rm d}x \\
&\qquad =\int_{\mu  }^{\infty }  
\frac{W\left(\Phi
\left(x\right)\right)^2}{\Phi \left(\mu \right)^2}\; e^{-\Phi\left(x\right)/\Phi
\left(\mu \right)}\; {\rm d}\Phi\left(x\right)\\
&\qquad=\int_{\Phi \left(\mu \right) }^{\infty }  
\frac{W\left(z\right)^2}{\Phi \left(\mu \right)^2 }\;e^{-z/\Phi
\left(\mu \right)}\; {\rm d}z =\int_{1}^{\infty
}w\left(v\right)^2\;e^{-v}\;{\rm d}v.
\end{align*}

\bigskip

To show the consistency of this test against any fixed alternative
 $$
{\scr H}_1:\quad S\left(x\right)\not=S_0\left(x\right)
$$ 
we just note that
$$
\sqrt{T}\left(\hat
f_T\left(x\right)-f_0\left(x\right)\right)=\sqrt{T}\left(\hat
f_T\left(x\right)-f\left(x\right)\right)+\sqrt{T}\left(
f\left(x\right)-f_0\left(x\right)\right) 
$$
where the first term is asymptotically normal and for the second term we have
$$
T\int_{\mu }^{\infty }h\left(x\right)\left(f\left(x\right)-f_0\left(x\right) \right)^2{\rm d}F_0\left(x\right)\longrightarrow \infty .
$$
Here $f\left(x\right)$ is the invariant density function (under
alternative). Of course, we have to suppose that the corresponding integrals
like \eqref{3} and \eqref{4} are finite. 

\bigskip

We now apply the similar arguments to study the Kolmogorov-Smirnov type statistic
$$
\gamma _T=\sup_{x\geq \mu }\sqrt{T}g\left(x\right)\left|\hat
f_T\left(x\right)-f_0\left(x\right) \right| .
$$
Our goal is to chose such weight function $g\left(\cdot \right)\geq 0$ that
the GoF test $\hat \psi _T=1_{\left\{\gamma _T>e_\varepsilon \right\}}$ based on
this statistic be ADF. Let us put $\sigma
\left(x\right)\equiv 1$ (for simplicity). Remind that the weak
convergence
$$
\zeta _T\left(\cdot \right)\Longrightarrow \zeta \left(\cdot \right)
$$
in the space of continuous functions vanishing at infinity was already proved
in \cite{Kut04}, Theorem 4.13. This convergence provides as well the
convergence of our statistic $\gamma _T$ with $g\left(x\right)\equiv
1$. Suppose that the function $g\left(\cdot \right)$ is such that we have the
convergence $\gamma _T\Rightarrow \gamma _0=\sup_xg\left(x \right)\left|\zeta
\left(x\right)\right|$ too. Let us put
$$
g\left(x\right)=\frac{1}{2f_0\left(x\right)\sqrt{\Phi \left(\mu
\right)}}\;e^{-\Phi \left(x\right)/\Phi \left(\mu \right)} .
$$
Then we can write (egalities in distribution)
\begin{align*}
 \gamma _0&=\sup_{x\geq \mu } 2g\left(x \right)f_0\left(x\right)\,\left|\int_{-\infty
}^{\infty } \frac{ 1_{\left\{y>x\right\}}-F_0\left(y\right)}{\sigma
\left(y\right){\sqrt{f_0\left(y\right)}}}\; {\rm d}W_y\right|\\
&=\sup_{x\geq \mu } 2g\left(x \right)f_0\left(x\right)\,\left|W\left(\Phi
\left(x\right)\right)\right| \\
&=\sup_{x\geq \mu } \frac{\left|W\left(\Phi
\left(x\right)\right)\right|}{\sqrt{\Psi \left(\mu \right)}}\;e^{-\Phi
\left(x\right)/\Phi \left(\mu \right)}=\sup_{v\geq 1
}\left|w\left(v\right)\right|e^{-v}. 
\end{align*}
We see that the test $\psi_T=1_{\left\{\gamma _T\geq e_\varepsilon \right\}}$ with
$e_\varepsilon $ from the equation
$$
\Pb\left\{\sup_{v\geq 1} \left|w\left(v\right)\right|e^{-v}\geq e_\varepsilon
\right\}=\varepsilon  
$$
is ADF and belongs to ${\cal K}_\varepsilon $.

\bigskip

{\bf Remark.} Note that these arguments do not work directly in the case of
{\it double sided alternatives}, i.e., if the function $S\left(x\right)$
changes under alternative for the values $x<\mu $ too. It can be shown that
the limit 
$$
\delta _T\Longrightarrow \int_{1}^{\infty
}\left[w_1\left(v\right)^2+w_2\left(v\right)^2\right]e^{-v}{\rm d}v 
$$
holds, but the Wiener processes $w_1\left(\cdot \right)$ and $w_2\left(\cdot
\right)$ are correlated and the correlation function depends on the model.

\section{GoF test based on EDF}

We study the GoF test $\Psi _T=1_{\left\{\Delta _T>c_\varepsilon \right\}}$ with
Cram\'er-von Mises type statistic
$$
\Delta _T=T\int_{\mu  }^{\infty }H\left(x\right)\left(\hat
F_T\left(x\right)-F_0\left(x\right)\right)^2 \,{\rm d}F_0\left(x\right)
$$
and our goal is to chose such weights $H\left(\cdot \right)\geq 0$ that this
statistic converges to the {\it distribution free} limit:
$$
\Delta _T\Longrightarrow \int_{1}^{\infty }w\left(v\right)^2e^{-v}\,{\rm d}v.
$$

This statistic can be written in the {\it empirical form} like \eqref{emp}
$$
\Delta _T^*=\int_{0  }^{T }H\left(X_t\right)\left(\hat
F_T\left(X_t\right)-F_0\left(X_t\right)\right)^2 \,{\rm d}t
$$
because $\Delta _T=\Delta _T^*+o\left(1\right)$ with the same explication as
above.

The properties of this statistic are quite close to that of $\delta _T$, that
is why we do not give here all details.

We suppose that the conditions ${\scr ES, RP}$ and $\Phi \left(x\right)<\infty
$ are fulfilled and  
\begin{align}
\label{9}
&\int_{\mu }^{\infty
}H\left(x\right)f_0\left(x\right)\Ex_0\left(\frac{F_0\left(\xi
\right)F_0\left(x\right)-F_0\left(\xi \wedge x\right) }{\sigma \left(\xi
\right)f_0\left(\xi\right)}\right)^2{\rm d}x<\infty \\ 
\label{10}
&\int_{\mu }^{\infty
}H\left(x\right)f_0\left(x\right)\Ex_0\left(\int_{\mu}^{\xi }\frac{F_0\left(y
\right)F_0\left(x\right)-F_0\left(y \wedge x\right) }{\sigma \left(y
\right)^2f_0\left(y\right)}{\rm d}y\right)^2{\rm d}x<\infty
\end{align}
The process $\eta _T\left(x\right)=\sqrt{T}\left(\hat
F_T\left(x\right)-F_0\left(x\right) \right)$ admits the presentations 
 (see \cite{Kut04},  p. 85)
\begin{align}
\eta _T\left(x\right)&=\frac{2}{\sqrt{T}}\int_{0}^{T}\frac{
F_0\left(X_t\right)F_0\left(x\right)-F_0\left(X_t\wedge x\right)}{\sigma
\left(X_t\right)f_0\left(X_t\right)}\; {\rm d}W_t\nonumber\\ 
&\label{11}\qquad
+\frac{2}{\sqrt{T}}\int_{X_0}^{X_T}\frac{ F_0\left(v\wedge
x\right)-F_0\left(v\right)F_0\left(x\right)}{\sigma
\left(v\right)^2f_0\left(v\right)}\; {\rm d}v.
\end{align}
 Note that if $\Phi \left(x\right)<\infty $, then (see \cite{Kut04}, Remark 1.64)
\begin{align*}
d_F\left(x\right)^2=4\Ex_0\left(\frac{
F_0\left(\xi \right)F_0\left(x\right)-F_0\left(\xi\wedge x\right)}{\sigma
\left(\xi\right)f_0\left(\xi\right)} \right)^2<\infty 
\end{align*}
too and 
$$
\sqrt{T}\left(\hat
F_T\left(x\right)-F_0\left(x\right) \right)\Longrightarrow {\cal
N}\left(0,d_F\left(x\right)^2 \right). 
$$ 
It can be shown that (under mild regularity conditions) 
\begin{align*}
\Delta _T\Longrightarrow \Delta _0=\int_{\mu }^{\infty
}H\left(x\right)f_0\left(x\right)\hat\eta \left(x\right)^2 {\rm d}x,
\end{align*}
where
$$
\hat\eta \left(x\right)=\left[F_0\left(x\right)-1\right]\int_{-\infty
}^{x}\frac{F_0\left(y\right)}{\sigma
\left(y\right)\sqrt{f_0\left(y\right)}}{\rm
d}W_y+F_0\left(x\right)\int_{ x}^{\infty}\frac{F_0\left(y\right)-1}{\sigma
\left(y\right)\sqrt{f_0\left(y\right)}}{\rm d}W_y.
$$
Hence
\begin{align*}
\Delta _0=\int_{\mu }^{\infty
}H\left(x\right)f_0\left(x\right)\left[F_0\left(x\right)-1\right]^2W
\left(\Psi \left(x\right)\right)^2 {\rm d}x, 
\end{align*}
where
\begin{align*}
\Psi \left(x\right)=\int_{-\infty}^{x}\frac{F_0\left(y\right)^2}{\sigma
\left(y\right)^2{f_0\left(y\right)}}{\rm d}y+F_0\left(x\right)^2\int_{
x}^{\infty} \left(\frac{F_0\left(y\right)-1}{F_0\left(x\right)-1} \right)^2
\frac{{\rm d}y}{\sigma  
\left(y\right)^2{f_0\left(y\right)}}.
\end{align*}
Further
\begin{align*}
\Delta _0&=\int_{\mu }^{\infty
}\frac{H\left(x\right)f_0\left(x\right)\left[F_0\left(x\right)-1\right]^2}{\Psi'
\left(x\right) }W \left(\Psi \left(x\right)\right)^2\, {\rm d} \Psi \left(x\right)\\
&=\Psi \left(0\right)^{-2}\int_{\mu }^{\infty} W \left(\Psi \left(x\right)\right)^2e^{-\Psi
\left(x\right)/\Psi \left(\mu\right)}\, {\rm d} \Psi \left(x\right)\\
&=\int_{1}^{\infty }w\left(v\right)^2e^{-v}\,{\rm d}v ,
\end{align*}
where we put $\Psi \left(x\right)=v\Psi \left(\mu\right) $,
$w\left(v\right)=\Psi \left(\mu\right)^{-1/2}W\left(v \Psi \left(\mu\right)
\right)$ and
$$
H\left(x\right)=\frac{\Psi' \left(x\right) }{\Psi
\left(\mu\right)^{2}f_0\left(x\right)\left[F_0\left(x\right)-1\right]^2}\,\,e^{-\Psi
\left(x\right)/\Psi \left(0\right)}.
$$
Of course, we suppose that the function $\Psi \left(x\right),x\geq \mu$ is
strictly monotone and $\Psi \left(\infty \right)=\infty $. 

\section{Examples}
 
Let us consider two examples. The first one is\\
{\it Example 1.} {\bf Ornstein-Uhlenbeck process.} Suppose that the observed
process under the basic hypothesis is 
$$
{\rm d}X_t=-a\left(X_t-b\right)\,{\rm d}t+\sigma \,{\rm d}W_t,\quad
X_0,\;0\leq t\leq T,
$$
where $a>0$. The invariant density is Gaussian $f_0\left(x\right)\sim {\cal
N}\left(b,\frac{\sigma ^2}{2a}\right)$ with median $\mu =b$. To check the
conditions \eqref{3} and \eqref{4} we estimate first the asymptotics of $\Psi
\left(x\right)$ as $x\rightarrow \infty $. We have
\begin{align*}
\Phi \left(x\right)&=\int_{-\infty }^{\infty } \frac{\left(
1_{\left\{y>x\right\}}-F_0\left(y\right)\right)^2}{\sigma
^2f_0\left(y\right)}\;  {\rm d}y=\int_{-\infty }^{x} \frac{
1}{\sigma
^2f_0\left(y\right)}\;  {\rm d}y\;\left(1+o\left(1\right)\right)\\
&=\frac{c}{x}\;\;e^{\frac{ax^2}{\sigma ^2}}\;\left(1+o\left(1\right)\right).
\end{align*}
Hence the both conditions are fulfilled and the test $\psi _T=1_{\left\{\delta
_T>d_\varepsilon \right\}}$ with
$$
\delta _T=T\int_{b}^{\infty }\frac{2F_0\left(x\right)-1 }{4\Phi \left(\mu
\right)^2\sigma 
\left(x\right)^2f_0\left(x\right)^3}\,\;e^{-\Phi \left(x\right)/\Phi
\left(\mu \right)}\;\left(\hat f_T\left(x\right)
- \frac{\sqrt{a}}{\sigma \sqrt{\pi }}e^{-\frac{a\left(x-b\right)^2}{\sigma
^2}}\right)^2{\rm d}x
$$
is ADF. 

It is easy to see that the conditions \eqref{9} and \eqref{10} are fulfilled
too and the test
 $\Psi _T=1_{\left\{\Delta
_T>d_\varepsilon \right\}}$ with
$$
\Delta _T=T\int_{b}^{\infty }\frac{\Psi' \left(x\right) }{4\Psi
\left(\mu\right)^{2}\left[F_0\left(x\right)-1\right]^2}\,\,e^{-\Psi
\left(x\right)/\Psi \left(0\right)}\;\left(\hat F_T\left(x\right)
- F_0\left(x\right)\right)^2 \,{\rm d}x
$$
is ADF. 

\bigskip

{\it Example 2.} {\bf Simple switching process.} Suppose that the observed process under
hypothesis is 
$$
{\rm d}X_t=-a\; \sgn\left(X_t-b\right)\,{\rm d}t+\sigma \,{\rm d}W_t,\quad
X_0,\;\; 0\leq t\leq T
$$
where $a>0$. The process is ergodic with invariant density 
$$
f_0\left(x\right)=\frac{a}{\sigma ^2}\;\exp\left\{-\frac{2a}{\sigma
^2}\left|x-b\right|\right\} 
$$
and median $\mu =b$. The function 
$$
\Psi \left(x\right)=\int_{-\infty }^{x} \frac{ 1}{\sigma
^2f_0\left(y\right)}\; {\rm
d}y\;\left(1+o\left(1\right)\right)=a^{-1}\,e^{\frac{2a}{\sigma
}x}\;\left(1+o\left(1\right)\right) 
$$
as $x\rightarrow \infty $ and the conditions \eqref{3}, \eqref{4} and
\eqref{9}, \eqref{10} are fulfilled. The direct calculation shows that $\Psi
\left(x\right)$ is strictly monotone function. Therefore the corresponding
tests
$\psi _T$ and $\Psi _T$ are ADF.

The limit distribution of the test statistic $\delta _T$ with
$h\left(x\right)\equiv 1$ were studied by Gassem \cite{Gas08}, who obtained
the Karhunen-Loeve expansion for the limit Gaussian process $\zeta
\left(\cdot \right)$. 

\section{Composite hypotheses}

Suppose that the observed diffusion process (under hypothesis ${\scr H}_0$) is
$$
{\rm d}X_t=S\left(\vartheta ,X_t\right)\,{\rm d}t+\sigma
\left(X_t\right)\,{\rm d}W_t,\quad X_0,\;\; 0\leq t\leq T
$$
where  $\vartheta $ is unknown parameter $\vartheta \in
\left(a,b\right)$. Therefore the basic hypothesis is composite. The test
statistic can be 
$$
\delta_T=T\int_{\mu }^{\infty }h\left(\hat\vartheta_T ,x\right) \left(\hat
f_T\left(x\right)-f_0(\hat\vartheta_T ,x)\right)^2{\rm d}F_0(\hat\vartheta_T ,x) ,
$$
where $\hat\vartheta _T$ is some consistent and asymptotically normal
estimator of $\vartheta $ and $h\left(\vartheta ,x\right)$ is the same
function as before with obvious modification, say,
$f_0\left(x\right)=f_0\left(\vartheta ,x\right)$. Unfortunately the test based
on this statistic is no more ADF because its
limit distribution depends on the distribution of estimator. To compensate
this contribution of estimator we can modify this statistic as follows (see,
e.g., Koul \cite{Kou} for similar transformation in time series). Suppose that
$\hat\vartheta_T$ is the MLE and the corresponding regularity conditions are
fulfilled (see \cite{Kut04}, Theorem 2.8). Then we have
$$
\sqrt{T}\left(\hat\vartheta _T-\vartheta \right)={\rm I}\left(\vartheta
\right)^{-1}T^{-1/2} \int_{0}^{T}\frac{\dot S\left(\vartheta ,X_t\right)}{\sigma
\left(X_t\right)} \,{\rm d}W_t  +o\left(1\right),
$$
where dot means derivation w.r.t. $\vartheta $. 

We want to substitute here the MLE, but in this case the stochastic integral
is not well defined, that is why we first rewrite this integral in the
following form (It\^o formula)
\begin{align*}
R_T \left(\vartheta \right)&   =    \int_{0}^{T}\frac{\dot
S\left(\vartheta,X_s\right)}{\sigma \left(X_s\right)} \; {\rm d}W_s
=\int_{0}^{T}\frac{\dot S\left(\vartheta 
,X_s\right)}{\sigma \left(X_s\right)^2}  \left[{\rm d}X_s-S\left(\vartheta
,X_s\right){\rm d}s \right] 
\\
&=\int_{x_0}^{X_T}\frac{\dot S\left(\vartheta,y\right) }{\sigma
\left(y\right)^2} \; {\rm d}y- \int_{0}^{T}\frac{\dot
S'\left(\vartheta,X_s\right)\sigma \left(X_s\right) -2\dot
S\left(\vartheta,X_s\right)\sigma' \left(X_x\right)}{2\sigma \left(X_s\right)
}\;{\rm d}s\\
&\quad - \int_{0}^{T}\frac{\dot S\left(\vartheta
,X_s\right)\,S\left(\vartheta ,X_s\right)}{\sigma \left(X_s\right)^2}\;  {\rm d}s.
\end{align*}
Here prim means derivation w.r.t. $x$.
The last expression contains no stochastic integral and we use it as definition of
 $R _T\left(\vartheta \right)$, where we can put  $\hat\vartheta _T$. 
 Now we can introduce the test statistic
\begin{align*}
&\hat\delta _T=T\int_{\mu }^{\infty }h\left(\hat\vartheta_T ,x\right) \left(\hat
f_T\left(x\right)-f_0(\hat\vartheta_T ,x)\right.\\
&\qquad \qquad \qquad \qquad\left.+\dot f_0(\hat\vartheta_T
,x)\;{\rm I}(\hat\vartheta_T
)\;T^{-1}\;R _T(\hat \vartheta_T )\right)^2\;{\rm d}F_0(\hat\vartheta_T ,x) .
\end{align*}
Note that
\begin{align*}
&\hat f_T\left(x\right)-f_0(\hat\vartheta_T ,x)=\hat
f_T\left(x\right)-f_0(\vartheta ,x)+
f_0\left(\vartheta ,x\right)-f_0(\hat\vartheta_T ,x)\\
&\qquad  =\hat
f_T\left(x\right)-f_0(\vartheta ,x)-\dot f_0(\hat\vartheta_T
,x)(\hat\vartheta_T - \vartheta )\left(1+o\left(1\right)\right).
\end{align*}
 Hence
\begin{align*}
\hat\delta _T&=T\int_{\mu }^{\infty }h\left(\vartheta ,x\right) \left(\hat
f_T\left(x\right)-f_0(\vartheta ,x)\right)^2{\rm d}F_0(\vartheta
,x)+o\left(1\right)\\
&\Longrightarrow \int_{1}^{\infty }w\left(v\right)^2e^{-v}\,{\rm d}v
\end{align*}
and the test $\hat\psi _T=1_{\left\{\hat\delta _T>d_\varepsilon \right\}}$ is
ADF. 

We supposed here that the median $\mu $ does not depend on $\vartheta $ (as in
Example 1 with $a=\vartheta $ and $\mu =b$).

\bigskip

In the case of Example 2 the situation is different. Suppose that $\vartheta $
is the shift parameter:
$$
{\rm d}X_t=-a\,\sgn\left(X_t-\vartheta \right)\,{\rm d}t+\sigma \,{\rm
d}W_t,\quad X_0,\;\ 0\leq t\leq T.
$$
Then we can use the statistic
$$
\hat\delta_T=T\int_{\hat\vartheta_T }^{\infty }h\left(\hat\vartheta_T ,x\right)
\left(\hat f_T\left(x\right)-f_0(\hat\vartheta_T ,x)\right)^2{\rm
d}F_0(\hat\vartheta_T ,x),
$$
and it can be shown that
$$
\hat\delta _T\Longrightarrow \int_{1}^{\infty }w\left(v\right)^2e^{-v}\,{\rm d}v.
$$
Indeed, the MLE $\hat \vartheta _T$ converges to $\vartheta $ with the rate
$T$ (and not $\sqrt{T}$) (see \cite{Kut04}, Theorem 3.26) and its contribution
to the limit distribution of $\delta _T$ is negligible.

Let us see what happens under alternative 
$$
{\scr H}_1\quad :\qquad S\left(\cdot \right)=S_*\left(\cdot \right),   \qquad
S_*\left(\cdot \right)\in {\cal F}_+
$$
where the set
$$
{\cal F}_+=\left\{S\left(\cdot 
\right): \inf_{\vartheta \in \Theta }\left\|\frac{S\left(\vartheta ,\cdot
\right)-S\left(\cdot \right)}{\sigma \left(\cdot \right)} \right\|_*>0\right\} 
$$
  and we suppose that the function $S_*\left(\cdot \right)$  satisfies
  the conditions ${\cal ES}$ and ${\cal RP}$. Therefore the invariant density
  is $f_{S_*}\left(\cdot \right)$. Here  the norm
$$
\left\|h\left(\cdot \right)\right\|_*^2=\int_{-\infty }^{\infty
}h\left(x\right)^2\,f_{S_*}\left(x\right)\,{\rm d}x. 
$$
Note that the MLE in this ``misspecified situation'' converges to the value
$\vartheta_* $ which minimizes the Kullback-Leibner distance
\begin{equation}
\label{15}
\left\|\frac{S\left(\vartheta_*  ,\cdot
\right)-S\left(\cdot \right)}{\sigma \left(\cdot \right)}
\right\|_*=\inf_{\vartheta \in\Theta }\left\|\frac{S\left(\vartheta ,\cdot 
\right)-S\left(\cdot \right)}{\sigma \left(\cdot \right)} \right\|_*
\end{equation}
(see \cite{Kut04}, Section 2.6.1 for details).

Hence, under regularity conditions we have
$$
\hat f_T\left(x\right)-f_0\left(\hat\vartheta _T,x\right)\longrightarrow
f_{S_*}\left(x\right) -f_0\left(\vartheta _*,x\right),\quad {\rm
I}\left(\hat\vartheta _T\right) \rightarrow {\rm I}\left(\vartheta _*\right)
$$
and $\dot f_0\left(\hat\vartheta _T,x\right)\rightarrow \dot
f_0\left(\vartheta _*,x\right) $. Further, it can be shown that
\begin{align*}
\frac{R_T\left(\hat\vartheta _T\right)}{T}&=\frac{1}{{T}}\int_{0}^{T}\frac{\dot
S\left(\hat\vartheta
_T,X_t\right)\left[S_*\left(X_t\right)-S\left(\hat\vartheta
_T,X_t\right)\right]}{\sigma \left(X_t\right)^2}\;{\rm d}t
\left(1+o\left(1\right)\right)\\
&=\frac{1}{{T}}\int_{0}^{T}\frac{\dot
S\left(\vartheta
_*,X_t\right)\left[S_*\left(X_t\right)-S\left(\vartheta
_*,X_t\right)\right]}{\sigma \left(X_t\right)^2}\;{\rm d}t
\left(1+o\left(1\right)\right)\\
&\longrightarrow \int_{-\infty }^{\infty }\frac{\dot
S\left(\vartheta
_*,x\right)\left[S_*\left(x\right)-S\left(\vartheta
_*,x\right)\right]}{\sigma \left(x\right)^2}\;f_{S_*}\left(x\right)\;  {\rm
d}x=R\left(\vartheta _*\right). 
\end{align*}
Therefore,
$$
\hat\delta _T\sim T\int_{-\infty }^{\infty }h\left(\vartheta _*,x\right)\left[f_{S_*}\left(x\right) -f_0\left(\vartheta _*,x\right)+\dot
f_0\left(\vartheta _*,x\right){\rm I}\left(\vartheta _*\right)R\left(\vartheta _*\right)
\right]^2{\rm d}x 
$$
and this test can be non consistent against alternatives $S_*\left(\cdot \right)$ such that
$$
f_{S_*}\left(x\right) -f_0\left(\vartheta _*,x\right)+\dot
f_0\left(\vartheta _*,x\right){\rm I}^{-1}\left(\vartheta
_*\right)R\left(\vartheta _*\right)=0 .
$$
Suppose that $\vartheta _*$ is an interior point of $\Theta $ and show that
the last equality is impossible. The value $\vartheta _*$ defined by the
equation \eqref{15} is the same time one of the solutions of the equation
$$
\int_{-\infty }^{\infty }\frac{\dot S\left(\vartheta
_*,x\right)\left[S_*\left(x\right)-S\left(\vartheta _*,x\right)\right]}{\sigma
\left(x\right)^2}\;f_{S_*}\left(x\right)\; {\rm d}x=0.
$$
Hence $R\left(\vartheta _*\right)=0 $.  The equality $f_{S_*}\left(x\right)
=f_0\left(\vartheta _*,x\right) $ implies
$$
\int_{0}^{x}\frac{S_*\left(y\right)}{\sigma \left(y\right)^2}\; {\rm
d}y=\int_{0}^{x}\frac{S\left(\vartheta _*,y\right)}{\sigma \left(y\right)^2}\;
{\rm d}y,
$$
which gives us $S_*\left(x\right)=S\left(\vartheta _*,x\right) $  for almost
all $x$ and the last equality contradicts the definition of
alternative. Therefore, $\hat\delta _T\rightarrow \infty $ and the test
$\hat\psi _T$ is consistent. 

\section{Discussion}

Note, that the similar problems of ADF GoF tests for stochastic differential equations
with ``small noise'' are considered in \cite{Kut08}.

The tests $\psi _T=1_{\left\{\delta _T>d_\varepsilon \right\}}$ and $\Psi
_T=1_{\left\{\Delta _T>c_\varepsilon \right\}}$ studied in this work are
consistent against any fixed alternative and it can be easily shown, that
these tests are {\it uniformly consistent} if the alternatives are separated from
hypothesis as follows
$$
{\scr H}_1\quad :\qquad S\left(\cdot \right)\in {\cal F}_r=\left\{S\left(\cdot
\right): \left\|f_S\left(\cdot \right)-f_0\left(\cdot \right) \right\|\geq r\right\}
$$
 with some $r>0$ for $\psi _T$  or
$$
{\scr H}_1\quad :\qquad S\left(\cdot \right)\in {\cal F}_q=\left\{S\left(\cdot
\right): \left\|F_S\left(\cdot \right)-F_0\left(\cdot \right) \right\|\geq q\right\}
$$
for $\Psi _T$  with some $q>0$. Here the norm
$$
\left\|m \left(\cdot \right) \right\|^2=\int_{-\infty }^{\infty
}m\left(x\right)^2{\rm d}F_0\left(x\right). 
$$
This means that
$$
\inf_{S\left(\cdot \right) \in {\cal F}_r}\Pb_S\left\{\delta _T>d_\varepsilon
\right\}\longrightarrow 1,\qquad \inf_{S\left(\cdot \right) \in {\cal
F}_q}\Pb_S\left\{\Delta _T>c_\varepsilon 
\right\}\longrightarrow 1  . 
$$

But if the alternative is defined by the Kullback-Leibner distance ($s>0$)
$$
{\scr H}_1\quad :\qquad S\left(\cdot \right)\in {\cal H}_s=\left\{S\left(\cdot
\right): \left\|\frac{S\left(\cdot \right)-S_0\left(\cdot \right)}{\sigma
\left(\cdot \right)} \right\|\geq s\right\}, 
$$
then the both tests are no more uniformly consistent. For example, the functions
$$
S_n\left(x\right)=S_0\left(x\right)+\alpha \sigma \left(x\right)^2
\cos\left(nx\right),\quad  n=1,2,\ldots 
$$ 
can belong to ${\cal H}_s $ but
$$
\left\|f_{S_n}\left(\cdot \right)-f_{0}\left(\cdot \right)\right\|\rightarrow
0\quad {\rm as}\qquad n\rightarrow \infty 
$$
and
$$
\inf_{S_n \in {\cal H}_s}\Pb_{S_n}\left\{\delta _T>d_\varepsilon
\right\}\longrightarrow \varepsilon \quad {\rm as}\qquad T\rightarrow \infty .
$$

For such alternatives it is better to use the
Chi-squared tests, which can be even asymptotically optimal in minimax sense.
The construction of such tests for signals in white Gaussian noise can be
found in Ermakov \cite{Er} and  Ingster and Suslina \cite{IS-03}. For inhomogeneous Poisson
processes see \cite{IK-07}. It is interesting to study such tests
in the case of ergodic diffusion processes too.

\end{document}